\documentclass[12pt]{article}

\usepackage{amsfonts,amsmath,amssymb,bbm}
\usepackage{graphicx}

\usepackage[hyperfootnotes=false]{hyperref}
\usepackage{setspace} 
\def\proof{\noindent{\bf Proof:}\hskip10pt}        
\def\QED{\hfill $\Box$}

\font\tenmath=msbm10 scaled 1200
\font\sevenmath=msbm7 scaled 1200
\font\Fivemath=msbm5 scaled 1200
\newfam\mathfam \textfont\mathfam=\tenmath
\scriptfont\mathfam=\sevenmath \scriptscriptfont\mathfam=\Fivemath

\def \\ { \cr }
\def\R{\mathbb{R}}
\def \1{1 \mkern -6mu 1} 
\def\N{\mathbb{N}}
\def\E{\mathbb{E}}
\def\P{\mathbb{P}}

\def\R{\mathbb{R}}

\def \e{{\rm e}}

\def \dt{{\ensuremath{\textup{d}}}}
\def \pB{{\ensuremath{\mathcal{B}}}}

\def \pU{{\ensuremath{\mathcal{U}}}}
\def \pY{{\ensuremath{\mathcal{Y}}}}
\def \pZ{{\ensuremath{\mathcal{Z}}}}

\newtheorem{theorem}{Theorem} 
\newtheorem{proposition}{Proposition} \newtheorem{lemma}{Lemma}
\newtheorem{corollary}{Corollary} 

\begin{document}

\title{Cutting edges at random in large recursive trees} \author{{Erich
    Baur\footnote{erich.baur@math.uzh.ch} { and } Jean
    Bertoin\footnote{jean.bertoin@math.uzh.ch}}\\ ENS Lyon and
  Universit\"at Z\"urich}
\maketitle 
\thispagestyle{empty}
\begin{abstract}

  We comment on old and new results related to the destruction of a random
  recursive tree (RRT), in which its edges are cut one after the other in a
  uniform random order. In particular, we study the number of steps
  needed to isolate or disconnect certain distinguished vertices
  when the size of the tree tends to infinity. New probabilistic explanations are
  given in terms of the so-called cut-tree and the tree of component sizes,
  which both encode different aspects of the destruction process. Finally,
  we establish the connection to Bernoulli bond percolation on large RRT's
  and present recent results on the cluster sizes in the supercritical
  regime.
\end{abstract}
\footnote{{\it Acknowledgment of support.} The research of the first author was
    supported by the Swiss National Science Foundation grant
    P2ZHP21\_51640.}
\section{Introduction}
\label{Sintro}
Imagine that we destroy a connected graph by removing or cutting its edges
one after the other, in a uniform random order.  The study of such a
procedure was initiated by Meir and Moon in \cite{MM1}. They were
interested in the number of steps needed to isolate a distinguished vertex
in a (random) Cayley tree, when the edges are removed uniformly at random from the
current component containing this vertex. Later on, Meir and Moon
\cite{MM2} extended their analysis to random recursive trees. The latter
form an important family of increasing labeled trees (see Section
\ref{Smaintools} for the definition), and it is the goal of this paper to
shed light on issues related to the destruction of such trees.

Mahmoud and Smythe~\cite{MaSmy} surveyed a multitude of results and
applications for random recursive trees. Their recursive structure make
them particularly amenable to mathematical analysis, from both a
combinatorial and probabilistic point of view. We focus on the
probabilistic side. Our main tools include the fundamental splitting
property, a coupling due to Iksanov and M\"ohle \cite{IM} and the so-called
{\it cut-tree} (see \cite{Be3}), which records the key information about
the destruction process. The cut-tree allows us to re-prove the results of
Kuba and Panholzer \cite{KP1} on the multiple isolation of
nodes. Moreover, we gain information on the number of steps needed to disconnect a
finite family nodes. 

Finally, we relate the destruction of a random recursive tree to Bernoulli
bond percolation on the same tree. We explain some results concerning the
sizes of percolation clusters in the supercritical regime, where the root
cluster forms the unique giant cluster.
\begin{figure}
\begin{center}\parbox{6cm}{\includegraphics[width=5cm]{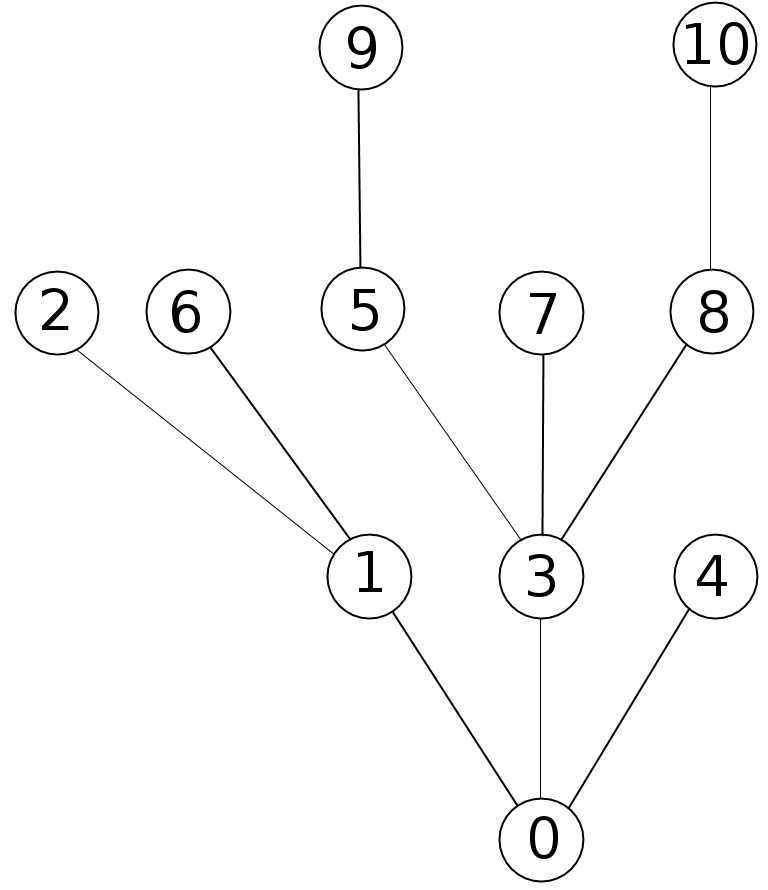}}
\end{center}
\centerline{\bf Figure 1}
\centering{\sl A recursive tree on the vertex set $\{0,1,\ldots,10\}$.}
\end{figure}

\section{Main tools}
\label{Smaintools}
In this section, we present some basic tools in the study of random
recursive trees  which will be useful to our purposes.

\subsection{The recursive construction, Yule process and P\'olya urn}

Consider a finite and totally ordered set of vertices, say $V$. A tree on
$V$ is naturally rooted at the smallest element of $V$, and is called
increasing if and only if the sequence of vertices along a segment from the
root to an arbitrary vertex increases. Most of the time we shall take
$V=\{0,1,\ldots, n\}$, which induces of course no loss of generality. More
precisely, it is convenient to introduce the following notion. For an
arbitrary totally ordered set $V$ with cardinality $|V|=n+1$, we call the
bijective map from $V$ to $\{0,1, \ldots,n\}$ which preserves the order,
the canonical relabeling of vertices.  Plainly the canonical relabeling
transforms an increasing tree on $V$ into an increasing tree on $\{0,1,
\ldots, n\}$. Such relabelings enable us to focus on the structure of the
rooted tree without retaining specifically the elements of $V$.

A {\it random recursive tree} (in short, RRT) on $\{0,1,\ldots, n\}$ is a
tree picked uniformly at random amongst all the increasing trees on
$\{0,1,\ldots, n\}$; it shall be denoted henceforth by $T_n$. In
particular, $T_n$ has $n$ edges and size (i.e. number of vertices)
$|T_n|=n+1$.  The terminology stems from the easy observation that a
version of $T_n$ can be constructed by the following simple recursive
random algorithm in which vertices are incorporated one after the
other. The vertex $1$ is naturally connected by an edge to the root $0$,
then $2$ is connected either to $0$ or to $1$ with equal probability $1/2$,
and more generally, the parent of the vertex $i$ is chosen uniformly at
random amongst $0,1, \ldots, i-1$ and independently of the other
vertices. This recursive construction is a close relative to the famous
Chinese Restaurant construction of uniform random permutations (see, for
instance, Section 3.1 in Pitman \cite{PiSF}), and in particular the number
of increasing trees of size $n+1$ equals $n!$.

Another useful observation is that this recursive construction can be
interpreted in terms of the genealogy of a Yule process. Recall that a Yule
process describes the evolution in continuous time of a pure birth process
in which each individual gives birth to a child at unit rate and
independently of the other individuals. We label individuals in the
increasing order of their birth times, the ancestor receiving by convention
the label $0$. If we let the process evolve until the population reaches
size $n+1$, then its genealogical tree, that is the tree where individuals
are viewed as vertices and edges connect children to their parent, is
clearly a RRT. Here is an application to percolation on $T_n$ which will be
useful later on.

\begin{lemma}\label{L2}
   Perform a Bernoulli bond percolation on $T_n$ with parameter $0<p<1$
   (i.e. each edge of $T_n$ is deleted with probability $1-p$,
   independently of the other edges), and let $C_n^0(p)$ denote the size of
   the cluster containing the root. Then
$$\lim_{n\to\infty} n^{-p}C_n^0(p) = C^0(p)\qquad \hbox{in distribution,}$$
where $C^0(p)>0$ a.s. is some random variable.
\end{lemma}

\proof We view $T_n$ as the genealogical tree of a standard Yule process
$(\pY_s)_{s\geq 0}$ up to time $\rho_n=\inf\{s\geq 0: \pY_s=n+1\}$. It is
well-known that the process $\e^{-s}\pY_s$ is a martingale which converges
a.s. to some random variable $W$ with the exponential distribution, and it
follows that
$$\lim_{n\to \infty} n^{-1}\e^{\rho_n}=1/W \qquad \hbox{a.s.}$$

In this setting, performing a Bernoulli bond percolation can be interpreted
as superposing neutral mutations to the genealogical tree, namely each
child is a clone of its parent with probability $p$ and a mutant with a new
genetic type with probability $1-p$, independently of the other
children. Neutrality means that the rate of birth does not depend on the
genetic type. Then the process $(\pY_s(p))_{s\geq 0}$ of the number of
individuals with the same genetic type as the ancestor is again a Yule
process, but now with birth rate $p$. As a consequence
$$\lim_{s\to \infty} \e^{-ps}\pY_s(p) = W(p)\qquad \hbox{a.s.},$$
where $W(p)$ denotes another exponentially distributed random variable.
We then observe that 
$$C_n^0(p)= \pY_{\rho_n}(p)\sim W(p)\e^{p\rho_n} \sim W(p)W^{-p}n^p,$$
which completes the proof. \QED

Plainly, the recursive construction can also be interpreted in terms of
urns, and we conclude this section by exemplifying this
connection. Specifically, the size of the root cluster $C_n^0(p)$ in the
above lemma can be identified as the number of red balls in the following
P\'olya-Hoppe urn. Start with one red ball which represents the root of the
tree. A draw is effected as follows:  \linebreak
(i) Choose a ball at random from the
urn, observe its color, and put the ball back to the urn. \linebreak
(ii) If its color
was red, add a red ball to the urn with probability $p$, and add a black
ball to the urn with probability $1-p$. If its color was black, add another
black ball to the urn. Then, after $n$ draws, the number of red balls is
given by $C_n^0(p)$, and in this way, Lemma \ref{L2} yields a limit theorem
for the proportion of red balls.

The choice $p=1$ in this urn scheme corresponds to the usual P\'olya
urn. Here, if one starts with one red ball and $k$ black balls, then the
number of red balls after $n-k$ draws is distributed as the size of the
subtree $T_n^k$ of a RRT $T_n$ that stems from the vertex $k$. It is
well-known from the theory of P\'olya urns that this number follows the
beta-binomial distribution with parameters $(n-k, 1, k)$. Moreover,
\begin{equation}
\label{EPolya}
\lim_{n\rightarrow\infty}n^{-1}|T_n^k|=\beta(1,k)\qquad \hbox{in distribution},
\end{equation}
where $\beta(1,k)$ is a beta$(1,k)$-distributed random variable. We will use
this fact several times below.

\subsection{The splitting property}
The {\it splitting property} (also called randomness preserving property)
reveals the fractal nature of RTT's: roughly speaking, if one removes an
edge from a RRT, then the two subtrees resulting from the split are in
turn, conditionally on their sizes, independent RRT's. This is of course of
crucial importance when investigating the destruction of a RRT, as we can
then apply iteratively the splitting property when removing the edges
uniformly at random and one after the other.

We select an edge of $T_n$ uniformly at random and remove it. Then $T_n$
splits into two subtrees, say $\tau^0_n$ and $\tau^*_n$, where $\tau^0_n$
contains the root $0$.  We denote by $T^0_{n}$ and $T^*_{n}$ the pair of
increasing trees which then result from the canonical relabelings of the
vertices of $\tau^0_n$ and $\tau^*_n$, respectively.  Introduce also an
integer-valued variable $\xi$ with distribution
\begin{equation}\label{E2}
  \P(\xi=j)=\frac{1}{j(j+1)},\qquad j=1,2, \ldots
\end{equation}

\begin{proposition} \label{P1} {\rm (Meir and Moon \cite{MM2})} In the
  notation above, $|\tau_n^*|=|T^*_n|$ has the same law as $\xi$
  conditioned on $\xi\leq n$, that is
 $$\P(|T^*_n|=j)=\frac{n+1}{nj(j+1)},\qquad j=1,2, \ldots, n.$$
 Further, conditionally on $|T^*_{n}|=j$, $T^0_{n}$ and $T^*_{n}$ are two
 independent RRT's with respective sizes $n-j+1$ and $j$.
\end{proposition}

\proof There are $n n!$ configurations $({\bf t}, e)$ given by an
increasing tree ${\bf t}$ on $\{0,1, \ldots, n\}$ and a distinguished edge
$e$.  We remove the edge $e$ and then relabel vertices canonically in each
of the resulting subtrees. Let us enumerate the configurations that yield a
given pair $({\bf t}^0, {\bf t}^*)$ of increasing trees on $\{0,1,\ldots,
n-j\}$ and $\{0,1,\ldots, j-1\}$, respectively.

Let $k\in\{0,1, \ldots, n-1\}$ denote the extremity of the edge $e$ which
is the closest to the root $0$ in ${\bf t}$, and $V^*$ the set of vertices
which are disconnected from $k$ when $e$ is removed. Since ${\bf t}$ is
increasing, all the vertices in $V^*$ must be larger than $k$, and since we
want $|V^*|=j$, there are $\left(^{n-k}_{\ j}\right)$ ways of choosing
$V^*$ (note that this is possible if and only if $k\leq n-j$).  There are a
unique increasing tree structure on $V^*$ and a unique increasing tree
structure on $\{0,1, \ldots, n\}\backslash V^*$ that yield respectively
${\bf t}^*$ and ${\bf t}^0$ after the canonical relabelings.
 
Conversely, given ${\bf t}^0$, ${\bf t}^*$, $k\in\{0,1, \ldots, n-j\}$ and
$V^*\subset\{k+1, \ldots, n\}$ with $|V^*|=j$, there is clearly a unique
configuration $({\bf t}, e)$ which yields the quadruple $(k,V^*,{\bf
  t}^0,{\bf t}^*)$. Namely, relabeling vertices in ${\bf t}^0$ and ${\bf
  t}^*$ produces two increasing tree structures $\tau^0$ and $\tau^*$ on
$\{0,1, \ldots, n\}\backslash V^*$ and $V^*$, respectively. We let $e$
denote the edge $(k,\min V^*)$ and then ${\bf t}$ is the increasing tree
obtained by connecting $\tau^0$ and $\tau^*$ using $e$.
 
It follows from the analysis above that
 $$\P(T^0_n={\bf t}^0, T^*_n={\bf t}^*)= \frac{1}{n n!}\sum_{k=0}^{n-j}\left(^{n-k}_{\ j}\right).$$
 Now recall that
 $$\sum_{k=0}^{n-j}\left(^{n-k}_{\ j}\right) = \sum_{\ell=j}^{n}\left(^{\ell}_{j}\right)= \left(^{n+1}_{j+1}\right)$$
 to conclude that
$$\P(T^0_n={\bf t}^0, T^*_n={\bf t}^*)= \frac{n+1}{n(n-j)! (j+1)!} = \frac{n+1}{nj(j+1)}  \times \frac{1}{(n-j)! (j-1)!}.$$
Since there are $(n-j)!$ increasing trees with size $n-j+1$ and $(j-1)!$
increasing trees with size $j$, this yields the claim. \QED

\noindent{\bf Remark.}  It can be easily checked that the splitting
property holds more generally when one removes a fixed edge, that is the
edge connecting a given vertex $k\in\{1,\ldots, n\}$ to its parent. Of
course, the distribution of the sizes of the resulting subtrees then
changes; see the connection to P\'olya urns mentioned in the beginning.

\subsection{The coupling of Iksanov and M\"ohle }

The splitting property was used by Meir and Moon \cite{MM2} to investigate
the following random algorithm for isolating the root $0$ of a
RRT. Starting from $T_n$, remove a first edge chosen uniformly at random
and discard the subtree which does not contain the root $0$. Iterate the
procedure with the subtree containing $0$ until the root is finally
isolated, and denote by ${X}_n$ the number of steps of this random
algorithm. In other words, $X_n$ is the number of random cuts that are
needed to isolate $0$ in $T_n$.

Iksanov and M\"ohle \cite{IM} derived from Proposition \ref{P1} a useful
coupling involving an increasing random walk with step distribution given
by \eqref{E2}.  Specifically, let $\xi_1, \xi_2, \ldots$ denote a sequence
of i.i.d. copies of $\xi$ and set $S_0=0$, 
\begin{equation}\label{ESn}
S_n=\xi_1+\cdots + \xi_n.
\end{equation} 
Further, introduce the last time that the random walk $S$ remains below the level $n$, 
\begin{equation}\label{ELn}
L(n)=\max\{k\geq 0: S_k\leq n\}.
\end{equation}

\begin{corollary} \label{C1} {\rm (Iksanov and M\"ohle \cite{IM})} One can
  construct on the same probability space a random recursive tree $T_n$
  together with the random algorithm of isolation of the root, and a
  version of the random walk $S$, such that if
\begin{equation}\label{Enest}
T^0_{n,0}=T_n \supset T^0_{n,1}\supset \cdots \supset T^0_{n,X_n}=\{0\}
\end{equation}
 denotes the nested sequence of the subtrees containing the root induced by the algorithm, then
 $X_n\geq L(n)$ and  
\begin{equation}\label{E6}
(|T^0_{n,0}\backslash T^0_{n,1}|,\ldots, |T^0_{n,L(n)-1} \backslash T^0_{n,L(n)}|) = (\xi_1,\ldots, \xi_{L(n)}).
\end{equation}
\end{corollary}

\proof Let us agree for convenience that $T^0_{n,j}=\{0\}$ for every
$j>X_n$, and first work conditionally on $(|T^{0}_{n,i}|)_{i\geq
  1}$. Introduce a sequence $((\varepsilon_i, \eta_i))_{i\geq 1}$ of
independent pairs of random variables such that for each $i$,
$\varepsilon_i$ has the Bernoulli law with parameter
$1/|T^0_{n,i-1}|=\P(\xi\geq |T^0_{n,i-1}|)$ and $\eta_i$ is an independent
variable distributed as $\xi$ conditioned on $\xi\geq |T^0_{n,i-1}|$. Then
define for every $i\geq 1$
$$\xi_i=\left\{ \begin{matrix} |T^0_{n,i-1}|-|T^0_{n,i}| &\hbox{ if }& \varepsilon_i=0\\
\eta_i&\hbox{ if }& \varepsilon_i=1 \end{matrix}\right.
$$
and the partial sums $S_i=\xi_1+\cdots+ \xi_i$. Observe that
$\varepsilon_i=1$ if and only if $\xi_i\geq |T^0_{n,i-1}|$, and hence, by
construction, there is the identity
$$\min\{i\geq 1: \varepsilon_i=1\}=\min\{i\geq 1: S_i\geq n+1\}.$$
Therefore, \eqref{E6} follows if we show that $\xi_1,\xi_2 \ldots$ are
(unconditionally) i.i.d. copies of $\xi$. This is essentially a consequence
of the splitting property. Specifically, for $j\leq n$, we have
$$\P(\xi_1=j)=\P(\varepsilon_1=0)\P(n+1-|T^0_{n,1}|=j)=\frac{n}{n+1}\P(|T^*_n|=j)= \frac{1}{j(j+1)},$$
where we used the notation and the result in Proposition \ref{P1}, whereas
for $j>n$ we have $$\P(\xi_1=j)=\P(\varepsilon_1=1)\P(\xi=j\mid \xi\geq
n+1)= \frac{1}{j(j+1)}.$$ Next, consider the conditional law of $\xi_2$
given $\xi_1$ and $|T^0_{n,1}|$. Of course, $|T^0_{n,1}|\geq n+1-\xi_1$,
and this inequality is in fact an equality whenever $\xi_1\leq n$. We know
from the splitting property that conditionally on its size, say
$|T^0_{n,1}|=m+1$ with $m\leq n-1$, $T^0_{n,1}$ is a RRT. Therefore Proposition \ref{P1}
yields again for $j\leq m$
\begin{eqnarray*}
  &&\P\left(\xi_2=j\mid  \xi_1\hbox{ and }  |T^0_{n,1}|=m+1\right)\\
  &=&\P\left(\varepsilon_2=0\mid \xi_1\hbox{ and }
    |T^0_{n,1}|=m+1\right)\P(m+1-|T^0_{n,2}|=j\mid \xi_1\hbox{ and }
  |T^0_{n,1}|=m+1)\\ 
  &=&\frac{m}{m+1}\P(|T^*_m|=j)\\
  &=& \frac{1}{j(j+1)}. 
\end{eqnarray*}
Similarly for $j>m$
\begin{eqnarray*}
  &&\P(\xi_2=j\mid \xi_1\hbox{ and } |T^0_{n,1}|=m+1)\\
  &=&\P(\varepsilon_2=1\mid \xi_1\hbox{ and }  |T^0_{n,1}|=m+1)\P(\xi=j\mid \xi\geq m+1)\\
  &=& \frac{1}{j(j+1)}.
\end{eqnarray*}
This shows that $\xi_2$ has the same distribution as $\xi$ and is
independent of $\xi_1$ and $|T^0_{n,1}|$. Iterating this argument, we get
that the $\xi_i$ form a sequence of i.i.d. copies of $\xi$, which completes
the proof. \QED

\section{The number of random cuts needed  to isolate the root}
\label{Sisoroot}
Recall the algorithm of isolation of the root which was introduced in the
preceding section, and that $X_n$ denotes its number of steps for $T_n$,
i.e. $X_n$ is the number of random cuts that are needed to isolate the root $0$
in $T_n$.  Meir and Moon \cite{MM2} used Proposition \ref{P1} to investigate
the first two moments of ${X}_n$ and showed that
\begin{equation}\label{EMM}
\lim_{n\to \infty} \frac{ \ln n}{n}  {X}_n =1\qquad \hbox{in probability.}
\end{equation}
The problem of specifying the fluctuations of ${X}_n$ was left open until
the work by Drmota {\it et al.}, who obtained the following remarkable
result.

\begin{theorem} \label{T1} {\rm (Drmota, Iksanov, M\"ohle and R\"osler \cite{DIMR})} As $n\to \infty$, 
$$\frac{\ln^2 n}{n}{X}_n - \ln n - \ln\ln n$$
converges in distribution to a completely asymmetric Cauchy variable $X$
with characteristic function
\begin{equation}\label{E3}
\E(\exp(i t X))=\exp\left(it\ln |t| -\frac{\pi}{2}|t|\right),\qquad t\in\R.
\end{equation}
\end{theorem}
In short, the starting point of the proof in \cite{DIMR} is the identity in
distribution
\begin{equation}\label{E1}
{X}_n \stackrel{\rm (d)}{=}  1 + {X}_{n-D_n},
\end{equation}
where $D_n$ is a random variable with the law of $\xi$ given $\xi\leq n$,
and $D_n$ is assumed to be independent of $X_1, \ldots, X_n$. More
precisely, \eqref{E1} derives immediately from the splitting property
(Proposition \ref{P1}). Drmota {\it et al.}  deduce from \eqref{E1} a PDE
for the generating function of the variables $X_n$, and then singularity
analysis provides the key tool for investigating the asymptotic behavior of
this generating function and elucidating the asymptotic behavior of $X_n$.

Iksanov and M\"ohle \cite{IM} developed an elegant probabilistic argument
which explains the unusual rescaling and the Cauchy limit law in Theorem
\ref{T1}. We shall now sketch this argument.
\noindent{\bf Sketch proof of Theorem \ref{T1}:}\hskip10pt
One starts observing that the
distribution in \eqref{E2} belongs to the domain of attraction of a
completely asymmetric Cauchy variable $X$ whose law is determined by \eqref{E3}, namely
\begin{equation}\label{E4}
\lim_{n\to \infty} \left( n^{-1}S_n - \ln n \right) = -X \qquad \hbox{in distribution}.
\end{equation} 
Then one deduces from \eqref{E4} that the asymptotic behavior of the last-passage time \eqref{ELn} is given by
\begin{equation}\label{E5}
\lim_{n\to \infty} \left( \frac{\ln^2n}{n} L(n) - \ln n - \ln \ln n  \right) = X \qquad \hbox{in distribution},
\end{equation}
see Proposition 2 in \cite{IM}. This limit theorem resembles of course
Theorem \ref{T1}, and the relation between the two is explained by the
coupling of the algorithm of isolation of the root and the random walk $S$
stated in Corollary \ref{C1}, as we shall now see.

Let the algorithm for isolating the root run for $L(n)$ steps. Then
the size of the remaining subtree that contains the root is
$n+1-S_{L(n)}$, and as a consequence, there are the bounds
$$L(n)\leq X_n \leq L(n) + n-S_{L(n)},$$
since at most $\ell-1$ edge removals are needed to isolate the root in any
tree of size $\ell$.  On the other hand, specializing a renewal theorem of
Erickson \cite{Er} for the increasing random walk $S$, one gets that
$$\lim_{n\to \infty} \ln (n-S_{L(n)})/\ln n = U \qquad \hbox{in distribution},$$
where $U$ is a uniform $[0,1]$ random variable. In particular
$$\lim_{n\to \infty}  \frac{\ln^2n}{n}(n-S_{L(n)}) = 0 \qquad \hbox{in probability.}$$
Thus Theorem \ref{T1} follows from \eqref{E5}. 
\QED

It should be noted that there exists a vertex version of the isolation
algorithm, where one chooses a vertex at random and destroys it together
with its descending subtree. The algorithm continues until the root is
chosen. Using an appropriate coupling with $X_n$, one readily shows that
the number  of random vertex removals $X_n^{(v)}$ needed to
destroy a RRT $T_n$ satisfies $(X_n-X_n^{(v)})=o(n/\ln^2n)$ in probability.
Henceforth, we concentrate on cutting edges.

\noindent{\bf Remark.}
Weak limit theorems for the number of cuts to isolate the root vertex have
also been obtained for other tree models, like conditioned Galton-Watson
trees including e.g. uniform Cayley trees and random binary trees
(Panholzer \cite{Pa1} and, in greater generality, Janson \cite{Ja2}),
deterministic complete binary trees (Janson \cite{Ja1}) and random split
trees (Holmgren \cite{Ho1,Ho2}). More generally, Addario-Berry {\it et al.}
\cite{ABBH} and Bertoin \cite{Be0} found the asymptotic limit distribution
for the number of cuts required to isolate a fixed number $\ell\geq 1$ of
vertices picked uniformly at random in a uniform Cayley tree. This result
was further extended by Bertoin and Miermont \cite{BM} to conditioned
Galton-Watson trees. We point to the remark after Corollary \ref{C3} for more
on this. Turning back to RRT's, recent generalizations of Theorem \ref{T1}
were found first by Kuba and Panholzer \cite{KP0, KP1} and then by Bertoin \cite{Be3},
some of which will be discussed in the reminder of this paper.

In \cite{KP1}, Kuba and Panholzer considered the situation when one wishes
to isolate the {\it first} $\ell$ vertices of a RRT $T_n$, $0,1, \ldots, \ell-1$,
where $\ell\geq 1$ is fixed.  In this direction, one modifies the algorithm
of isolation of the sole root in an obvious way. A first edge picked
uniformly at random in $T_n$ is removed. If one of the two resulting
subtrees contains none of the vertices $0,1, \ldots, \ell-1$, then it is
discarded forever. Else, the two subtrees are kept. In both cases, one
iterates until each and every vertex $0,1, \ldots, \ell-1$ has been
isolated, and we write $X_{n,\ell}$ for the number of steps of this
algorithm.

The approach of Kuba and Panholzer follows analytic methods similar to the
original proof of Theorem \ref{T1} by Drmota {\it et al.} \cite{DIMR}.  We
point out here that the asymptotic behavior of $X_{n,\ell}$ can also be
deduced from Theorem \ref{T1} by a probabilistic argument based on the
following elementary observation, which enables us to couple the variables
$X_{n,\ell}$ for different values of $\ell$. Specifically, we run the usual
algorithm of isolation of the root, except that now, at each time when a
subtree becomes disconnected from the root, we keep it aside whenever it
contains at least one of the vertices $1, \ldots, \ell-1$, and discard it
forever otherwise. Once the root $0$ of $T_n$ has been isolated, we resume
with the subtree containing $1$ which was set aside, meaning that we run a
further algorithm on that subtree until its root $1$ has been isolated,
keeping aside the further subtrees disconnected from $1$ which contain at
least one of the vertices $2, \ldots, \ell-1$. We then continue with the
subtree containing the vertex $2$, and so on until each and every vertex
$0,1, \ldots, \ell-1$ has been isolated. If we write $X'_{n,\ell}$ for the
number of steps of this algorithm, then it should be plain that
$X'_{n,\ell}$ has the same law as $X_{n,\ell}$, and further
$X_n=X'_{n,1}\leq \cdots\leq X'_{n,\ell}$.

We shall now investigate the asymptotic behavior of the increments
$\Delta_{n,i}= X'_{n,i+1}- X'_{n,i}$ for $i\geq 1$ fixed.  In this
direction, suppose that we now remove the edges of $T_n$ one after the
other in a uniform random order until the edge connecting the vertex $i$ to
its parent is removed. Let $\tau^{i}_n$ denote the subtree containing $i$
that arises at this step.

\begin{lemma} \label{L1} For each fixed $i\geq 1$, 
$$\lim_{n\to\infty}\frac{\ln|\tau_n^{i}|}{\ln n}=U\qquad \hbox{in distribution},$$
where $U$ is a uniform $[0,1]$ random variable.
\end{lemma}

For the moment, let us take Lemma \ref{L1} for granted and deduce the following.

\begin{corollary} \label{C2} We have
$$\lim_{n\to\infty}\frac{\ln \Delta_{n,i}}{\ln n}=U\qquad \hbox{in distribution},$$
 where $U$ is a uniform $[0,1]$ random variable.
\end{corollary} 
\proof Just observe that $\Delta_{n,i}$ has the same law as the number of
cuts needed to isolate the root $i$ of $\tau_n^i$, and recall from an
iteration of the splitting property that conditionally on its size,
$\tau_n^i$ is a RRT. Our statement now follows readily from \eqref{EMM} and
Lemma \ref{L1}. \QED

Writing $X'_{n,\ell} = X_n + \Delta_{n,1}+\cdots + \Delta_{n,\ell-1}$, we
now see from Theorem \ref{T1} and Corollary \ref{C2} that for each fixed
$\ell\geq 1$, there is the weak convergence
\begin{equation}
\label{E8}
\lim_{n\to\infty}\left(\frac{\ln^2 n}{n}X'_{n,\ell} - \ln n - \ln\ln n\right) = X\qquad \hbox{in distribution,}
\end{equation}
which is Theorem 1 in \cite{KP1}. We now proceed to the proof of Lemma \ref{L1}.

\proof Let $T_n^i$ denote the subtree of $T_n$ that stems from the vertex
$i$, and equip each edge $e$ of $T_n$ with a uniform $[0,1]$ random
variable $U_e$, independently of the other edges. Imagine that the edge $e$
is removed at time $U_e$, and for every time $0\leq s \leq 1$, write
$T_n^i(s)$ for the subtree of $T_n^i$ which contains $i$ at time
$s$. Hence, if we write $U=U_e$ for $e$ the edge connecting $i$ to its
parent, then $ \tau^i_n=T_n^i(U)$.  Further, since $U$ is independent of
the other uniform variables, conditionally on $U$ and $T_n^i$, $\tau^i_n$
can be viewed as the cluster that contains the root vertex $i$ after a
Bernoulli bond percolation on $T^i_n$ with parameter $1-U$. Thus,
conditionally on $|T_n^i|=m+1$ and $U=1-p$, $|\tau^i_n|$ has the same law as
$C^0_m(p)$ in the notation of Lemma \ref{L2}.

From \eqref{EPolya} we know that $n^{-1}|T_n^i|$ converges in
distribution as $n\to \infty$ to a beta variable with parameters $(1,i)$,
say $\beta$, which is of course independent of $U$. On the other hand,
conditionally on its size, and after the usual canonical relabeling of its
vertices, $T_n^i$ is also a RRT (see the remark at the end of Section
\ref{Smaintools}). It then follows from Lemma \ref{L2} that
$$\lim_{n\to \infty}\frac{\ln |\tau_n^i|}{\ln n}=1-U\qquad \hbox{in probability,}$$
which establishes our claim. \QED

\section{The destruction process and its tree representations}
\label{Sdestructproc}
Imagine now that we remove the edges of $T_n$ one after the other and in a
uniform random order, no matter whether they belong to the root component
or not. We call this the {\it destruction process} of $T_n$. After $n$
steps, no edges are present anymore and all the vertices have been
isolated. In particular, the random variable which counts only the number
of edge removals from the root component can be identified
with $X_n$ from the previous section.

The purpose of this section is to introduce and study the asymptotic
behavior of two trees which can be naturally associated to this destruction
process, namely the tree of component sizes and the cut-tree. Furthermore, we
give some applications of the cut-tree to the isolation and disconnection
of nodes and comment on ordered destruction of a RRT.

\subsection{The tree of component sizes}
In this part, we are interested in the sizes of the tree
components produced by the destruction process. Our analysis
will also prove helpful for studying percolation clusters of a RRT in 
Section \ref{Spercolation}. 

The component sizes are naturally stored in a tree structure. As our index
set, we use the universal tree
\begin{equation*}
\pU = \bigcup_{k=0}^\infty \N^k,
\end{equation*}
with the convention $\N^0 = \{\emptyset\}$ and $\N=\{1,2,\ldots\}$. In
particular, an element $u\in\pU$ is a finite sequence of strictly positive
integers $(u_1,\ldots,u_k)$, and its length $|u|=k$ represents the
``generation'' of $u$. The $j$th child of $u$ is given by
$uj=(u_1,\ldots,u_k,j)$, $j\in\mathbb{N}$. The empty sequence $\emptyset$ is
the root of the tree and has length $|\emptyset|=0$. If no confusion
occurs, we drop the separating commas and write $(u_1\cdots u_k)$ or simply
$u_1\cdots u_k$ instead of $(u_1,\ldots,u_k)$. Also, 
$\emptyset u$ represents the element $u$.

We define a tree-indexed process $\pB^{(n)}=(\pB^{(n)}_u : u\in\pU)$, which
encodes the sizes of the tree components stemming from the destruction of
$T_n$. We will directly identify a vertex $u$ with its label
$\pB^{(n)}_u$. Following the steps of the destruction process, we build
this process dynamically starting from the singleton
$\pB^{(n)}_\emptyset=n+1$ and ending after $n$ steps with the full process
$\pB^{(n)}.$ More precisely, when the first edge of $T_n$ is removed in the
destruction process, $T_n$ splits into two subtrees, say $\tau_n^0$ and
$\tau_n^\ast$, where $\tau_n^0$ contains the root $0$.  We stress that
$\tau_n^0$ is naturally rooted at $0$ and $\tau_n^\ast$ at its smallest
vertex.  The size $|\tau_n^\ast|$ is viewed as the first child of
$\pB^{(n)}_\emptyset$ and denoted by $\pB_1^{(n)}$. Now first suppose that
the next edge which is removed connects two vertices in
$\tau_n^\ast$. Then, $\tau_n^\ast$ splits into two tree components. The
size of the component not containing the root of $\tau_n^\ast$ is viewed as
the first child of $\pB^{(n)}_1$ and denoted by $\pB^{(n)}_{11}$. On the
other hand, if the second edge which is removed connects two vertices in
$\tau_n^0$, then the size of the component not containing $0$ is viewed as
the second child of $\pB^{(n)}_\emptyset$ and denoted by $\pB_2^{(n)}$. It
should now be plain how to iterate this construction. After $n$ steps, we
have in this way defined $n+1$ variables $\pB_u^{(n)}$ with $|u|\leq n$,
and we extend the definition to the full universal tree by letting
$\pB_u^{(n)}=0$ for all the remaining $u\in\pU$. We refer to Figure $2$ for
an example. The tree components whose sizes are encoded by the elements
$\pB_u^{(n)}$ with $|u|=k$ are called the components of generation $k$.

To sum up, every time an edge is removed in the destruction process, a
tree component $\tau_n$ splits into two subtrees, and we adjoin the size of
the subtree which does not contain the root of $\tau_n$ as a new child to
the vertex representing $\tau_n$. Note that the root $\pB_\emptyset^{(n)}$
has $X_n$ many nontrivial children, and they represent the sizes of the tree
components which were cut from the root one after the other in the
algorithm for isolating the root. 

\begin{figure}
\begin{center}\parbox{6cm}{\includegraphics[width=5cm]{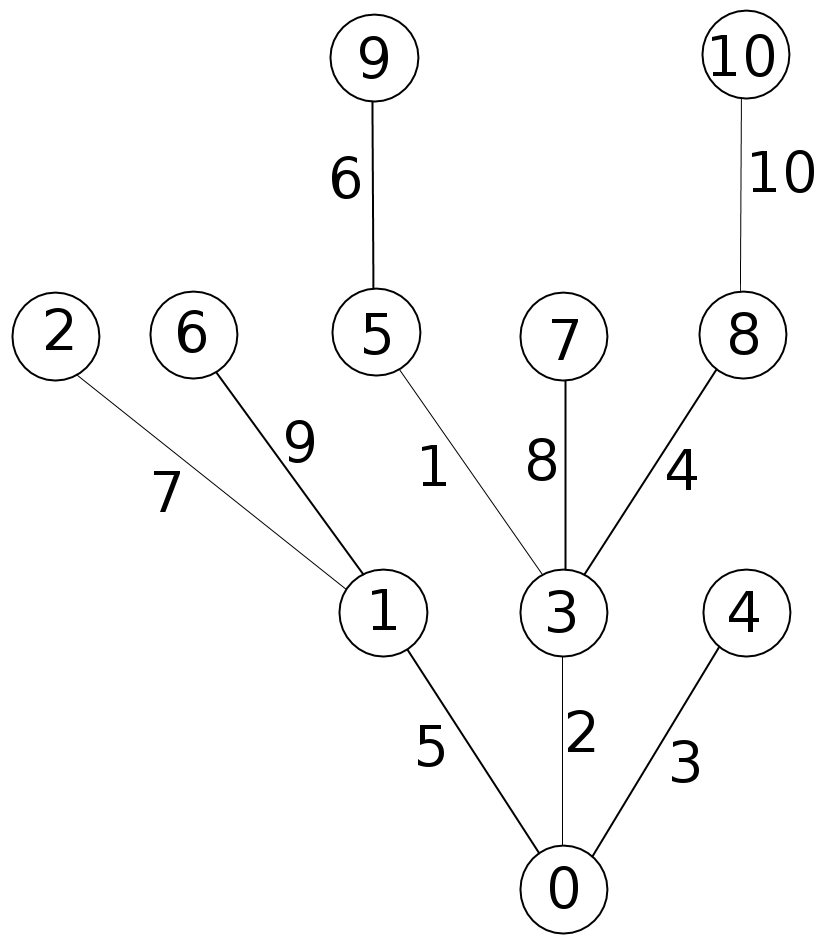}}
\parbox{9cm}{\includegraphics[width=9cm]{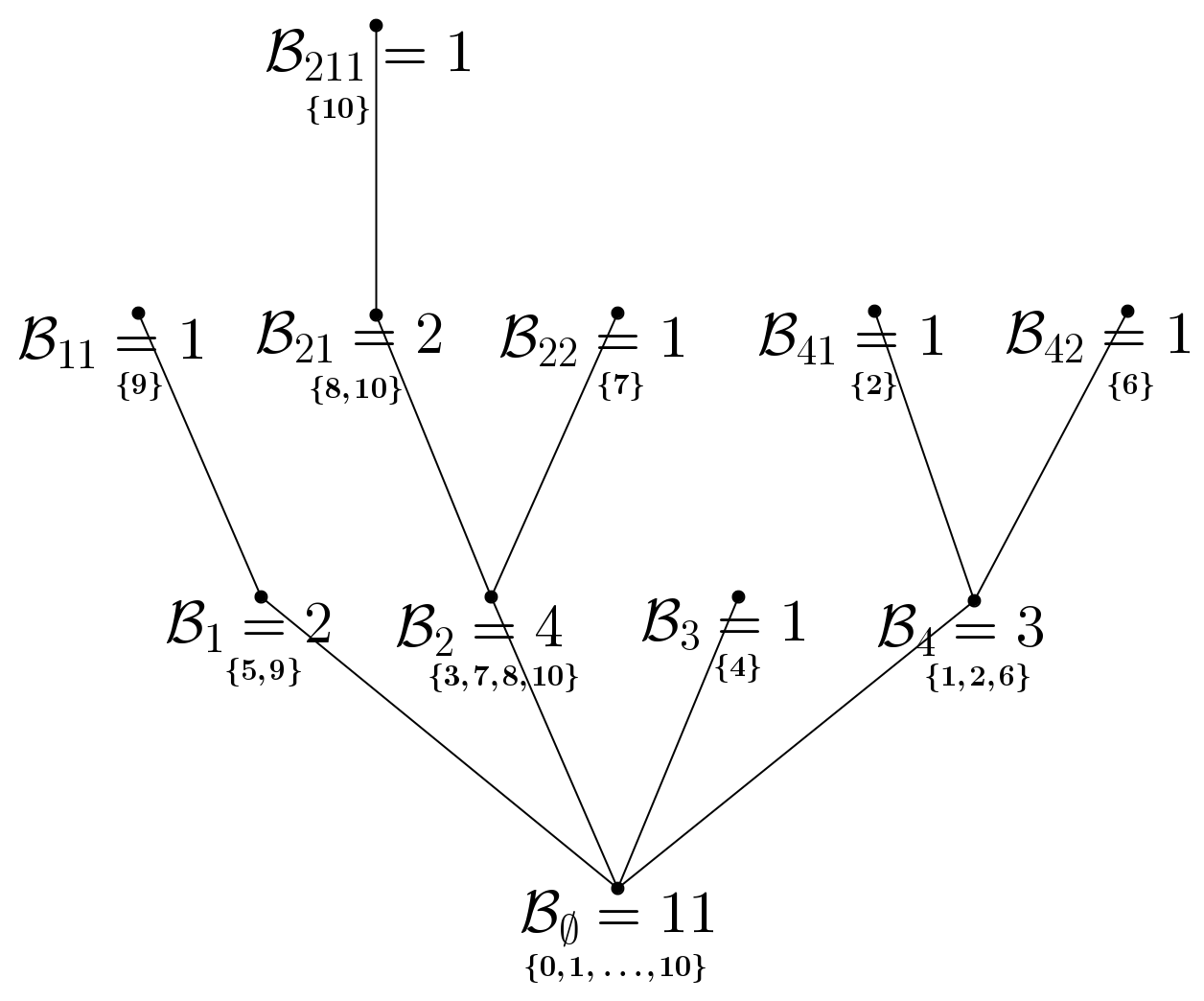}}
\end{center}
\centerline{\bf Figure 2}
\centering{\sl Left: A recursive tree with vertices labeled
  $0,1,\ldots,10$. The labels on the edges indicate the order in which they
  are removed by the destruction process.}\\
{\sl  Right: The corresponding tree of component sizes, with the vertex
  sets of the tree components. The elements $\pB_u^{(n)}$ of size $0$ are omitted.}
\end{figure}

We now interpret $\pB^{(n)}$ as the genealogical tree of a multi-type
population model, where the type reflects the size of the tree component
(and thus takes integer values). In particular the ancestor $\emptyset$ has
type $n+1$; furthermore, a node $u$ with $\pB^{(n)}_u=0$ corresponds to an
empty component and is therefore absent in the population model. We also
stress that the type of an individual is always given by the sum of the
types of its children plus $1$. As a consequence, types can be recovered
from the sole structure of the genealogical tree. More precisely, the type
of an individual is simply given by the total size of the subtree of the
genealogical tree stemming from that individual.

The splitting property of a RRT immediately transfers into a branching
property for this population model.

\begin{lemma}
\label{L3}  
The population model induced by the tree of component sizes $\pB^{(n)}$ is
a multi-type Galton-Watson process starting from one particle of type
$n+1$. The reproduction distribution $\lambda_i$ of an individual of type
$i\geq 1$ is given by the law of the sequence of the sizes of the non-root
subtrees which are produced in the algorithm for isolating the root of a
RRT of size $i$.
\end{lemma}

Even though the coupling of Iksanov and M\"ohle is not sufficient to fully
describe the reproduction law, it nonetheless provides essential
information on $\lambda_i$ in terms of a sequence of i.i.d. copies of
$\xi$. As we will see next, extreme value theory for the i.i.d. sequence
then enables us to specify asymptotics of the population model when the
type $n+1$ of the ancestor goes to infinity.

To give a precise statement, we rank the children of each
individual in the decreasing order of their types. Formally, given
the individual indexed by $u\in\pU$ has exactly $\ell$ children of type $\geq 1$, we
let $\sigma_u$ be the random permutation of $\{1,\ldots,\ell\}$ which sorts
the sequence of types $\pB^{(n)}_{u1},\ldots,\pB^{(n)}_{u\ell}$ in the decreasing
order, i.e.
\begin{equation*}
  \pB^{(n)}_{u\sigma_u(1)}\geq \pB^{(n)}_{u\sigma_u(2)}\geq \ldots\geq\pB^{(n)}_{u\sigma_u(\ell)},
\end{equation*}
where in the case of ties, children of the same type are ranked uniformly
at random. We extend $\sigma_u$ to a bijection $\sigma_u:\N\rightarrow\N$ by putting
$\sigma_u(i) = i$ for $i>\ell$.

We then define the
global random bijection $\sigma = \sigma^{(n)} : \pU\rightarrow\pU$
recursively by setting $\sigma(\emptyset) = \emptyset$,
$\sigma(j)=\sigma_\emptyset(j)$, and then, given $\sigma(u)$, $\sigma(uj) =
\sigma(u)\sigma_{\sigma(u)}(j)$, $u\in\pU$, $j\in\mathbb{N}$.  Note that
$\sigma$ preserves the parent-child relationship, i.e. children of $u$ are
mapped into children of $\sigma(u)$. We simply write $(\pB^{(n)\downarrow}_{u} :
u\in\pU)=(\pB^{(n)}_{\sigma(u)} : u\in\pU)$ for the process
which is ranked in this way.

Now, if the sizes of the components of generation $k$ are normalized by a factor $\ln^k n/n$,
we obtain finite-dimensional convergence of $\pB^{(n)\downarrow}$ towards
the genealogical tree of a continuous-state branching process with
reproduction measure $\nu(da)=a^{-2}\dt a$ on $(0,\infty)$. More precisely,
the limit object is a tree-indexed process $\pZ=(\pZ_u:u\in \pU)$ with
initial state $\pZ_\emptyset = 1$, whose distribution is characterized by
induction on the generations as follows.
\begin{enumerate}
\item $\pZ_\emptyset = 1$ almost surely;
\item for every $k=0,1,2,\ldots,$ conditionally on $(\pZ_v: v\in\pU,|v|\leq
  k)$, the sequences $(\pZ_{uj})_{j\in\N}$ for the vertices $u\in\pU$ at
  generation $|u|=k$ are independent, and each sequence
  $(\pZ_{uj})_{j\in\N}$ is distributed as the family of the atoms of a Poisson
  random measure on $(0,\infty)$ with intensity $\pZ_u\nu$, where the atoms
  are ranked in the decreasing order of their sizes.
 \end{enumerate}

\begin{proposition}
\label{Pbranching}
As $n\rightarrow\infty$, there is the convergence in the sense of
finite-dimensional distributions,
$$
\pZ^{(n)}=\left(\frac{(\ln n)^{|u|}}{n}\pB_{u}^{(n)\downarrow} :
  u\in\pU\right)\Longrightarrow \pZ.
$$
\end{proposition}
We only sketch the proof and refer to the forthcoming paper~\cite{Ba} for
details. Basically, if $\xi_1,\xi_2,\ldots$ is a sequence of of
i.i.d. copies of $\xi$, then for $a>0$, the number of indices $j\leq k$
such that $\xi_j\geq an/\ln n$ is binomially distributed with
parameters $k$ and $\lceil an/\ln n\rceil^{-1}$. From \eqref{E5} and
Theorem 16.16 in \cite{Ka} we deduce that for a fixed integer $j$, the $j$
largest among $\xi_1,\ldots,\xi_{L(n)}$, normalized by a factor $\ln n/n$,
converge in distribution to the $j$ largest atoms of a Poisson random
measure on $(0,\infty)$ with intensity $\nu(da)=a^{-2}\dt a$. Since
$n-S_{L(n)} = o(n/\ln^2 n)$ in probability, finite-dimensional convergence
of $\pZ^{(n)}$ restricted to generations $\leq 1$ then follows from
\eqref{E6}. Lemma \ref{L3} enables us to transport the arguments to the
next generations.

\subsection{The cut-tree}
\label{Scuttree}
Consider for a while a deterministic setting where $T$ is an arbitrary tree
on some finite set of vertices $V$. Imagine that its edges are removed one
after the other in some given order, so at the end of the process, all the
vertices of $T$ have been disconnected from each other. We shall encode the
destruction of $T$ by a rooted binary tree, which we call the cut-tree and
denote by ${\rm Cut}(T)$. The cut-tree has internal nodes given by the
non-singleton connected components which arise during the destruction, and
leaves which correspond to the singletons and which can thus be identified
with the vertices in $V$. More precisely, the root of ${\rm Cut}(T)$ is
given by $V$, and when the first edge of $T$ is removed, disconnecting $V$
into, say, $V_1$ and $V_2$, then $V_1$ and $V_2$ are viewed as the two
children of $V$ and thus connected to $V$ by a pair of edges. Suppose that
the next edge which is removed connects two vertices in $V_1$, so removing
this second edge disconnects $V_1$ into, say $V_{1,1}$ and $V_{1,2}$. Then
$V_{1,1}$ and $V_{1,2}$ are viewed in turn as the two children of $V_1$. We
iterate in an obvious way, see Figure $3$ for an example.\footnote{For the
  sake of simplicity, this notation does not record the order in which the
  edges are removed, although the latter is of course crucial in the
  definition of the cut-tree. In this part, we are concerned with uniform
  random edge removal, while in the last part of this section, we look at
  ordered destruction of a RRT, where edges are removed in the order of
  their endpoints most distant from the root.}

\begin{figure}
\begin{picture}(400,280)(-10,-30)

\put(80,10){\circle{20}}
\put(40,10){\circle{20}}
\put(40,50){\circle{20}}
\put(120,50){\circle{20}}
\put(80,90){\circle{20}}
\put(160,90){\circle{20}}
\put(40,130){\circle{20}}
\put(120,130){\circle{20}}
\put(80,170){\circle{20}}

\put (40,20){\line(0,1){20}}
\put (50,10){\line(1,0){20}}
\put (85,18){\line(1,1){26}}
 
\put (115,58){\line(-1,1){26}}
\put (75,98){\line(-1,1){26}}
\put (85,98){\line(1,1){26}}
\put (125,58){\line(1,1){26}}
\put (80,100){\line(0,1){60}}

\put(80,10)
{\makebox(0,0){a}}
\put(40,10)
{\makebox(0,0){h}}
\put(120,50)
{\makebox(0,0){g}}
\put(40,50)
{\makebox(0,0){f}}
\put(80,90)
{\makebox(0,0){i}}
\put(160,90)
{\makebox(0,0){d}}
\put(40,130)
{\makebox(0,0){e}}
\put(120,130)
{\makebox(0,0){c}}
\put(80,170)
{\makebox(0,0){b}}

\put(103,29){\makebox(0,0){1}}
\put(60,4){\makebox(0,0){6}}
\put(34,30){\makebox(0,0){8}}
\put(100,63){\makebox(0,0){2}}
\put(102,105){\makebox(0,0){5}}
\put(137,65){\makebox(0,0){4}}
\put(55,110){\makebox(0,0){3}}
\put(75,130){\makebox(0,0){7}}


\put(290,00){\framebox(60,15){abcdefghi} }

\put(250,50){\framebox(20,15){afh} }
\put(350,50){\framebox(40,15){bcdegi} }

\put(230,90){\framebox(10,15){a} }
\put(260,90){\framebox(15,15){fh} }  
\put(230,130){\framebox(10,15){f} }
  \put(260,130){\framebox(10,15){h} }

\put(330,90){\framebox(30,15){bcei} }
\put(400,90){\framebox(15,15){dg} }

\put(300,130){\framebox(20,15){bci} }
\put(340,130){\framebox(10,15){e} }
\put(370,130){\framebox(10,15){g} }
\put(400,130){\framebox(10,15){d} }

\put(300,170){\framebox(15,15){bi} }
\put(340,170){\framebox(10,15){c} }

\put(270,210){\framebox(10,15){b} }
\put(310,210){\framebox(10,15){i} }

\put (300,16){\line(-1,1){34}}
\put (330,16){\line(1,1){34}}

\put (260,66){\line(-1,1){24}}
\put (265,66){\line(0,1){24}}

\put (260,106){\line(-1,1){24}}
\put (265,106){\line(0,1){24}}

\put (365,66){\line(-1,1){24}}
\put (383,66){\line(1,1){24}}

\put (335,106){\line(-1,1){24}}
\put (345,106){\line(0,1){24}}

\put (400,106){\line(-1,1){24}}
\put (405,106){\line(0,1){24}}

\put (305,146){\line(0,1){24}}
\put (315,146){\line(1,1){24}}

\put (300,186){\line(-1,1){24}}
\put (313,186){\line(0,1){24}}

\end{picture}

\centerline{\bf Figure 3}
\centerline{\sl  Left: Tree $T$ with vertices labeled a,...,i; edges are
  enumerated in the order of the cuts.}
\centerline{\sl  Right: Cut-tree ${\rm Cut}(T)$ on the set of blocks
recording the destruction of  $T$.}
\end{figure}

It should be clear that the number of cuts required to isolate a given
vertex $v$ in the destruction of $T$ (as previously, we only count the cuts
occurring in the component which contains $v$) corresponds precisely to the
height of the leaf $\{v\}$ in ${\rm Cut}(T)$. More generally, the number of
cuts required to isolate $k$ distinct vertices $v_1, \ldots, v_k$ coincides
with the total length of the cut-tree reduced to its root and the $k$
leaves $\{v_1\}, \ldots, \{v_k\}$ minus $(k-1)$, where the length is
measured as usual by the graph distance on ${\rm Cut}(T)$. In short, the
cut-tree encapsulates all the information about the numbers of cuts needed
to isolate any subset of vertices.

We now return to our usual setting, that is $T_n$ is a RRT of size $n+1$,
whose edges are removed in a uniform random order, and we write ${\rm
  Cut}(T_n)$ for the corresponding cut-tree. We point out that the
genealogical tree of component sizes which was considered in the previous
section can easily be recovered from ${\rm Cut}(T_n)$. Specifically, the
root $\{0,1,\ldots, n\}$ of ${\rm Cut}(T_n)$ has to be viewed as the
ancestor of the population model, its type is of course $n+1$. Then the
blocks of ${\rm Cut}(T_n)$ which are connected by an edge to the segment
from the root $\{0,1,\ldots, n\}$ to the leaf $\{0\}$ are the children of
the ancestor in the population model, the type of a child being given by
the size of the corresponding block. The next generations of the population
model are then described similarly by an obvious iteration.

The segment of ${\rm Cut}(T_n)$ from its root $\{0,1,\ldots, n\}$ to the
leaf $\{0\}$ is described by the nested sequence \eqref{Enest}, and the
coupling of Iksanov and M\"ohle stated in Corollary \ref{C1} expresses the
sequence of the block-sizes along the portion of this segment starting from
the root and with length $L(n)$, in terms of the random walk $S$.  We shall
refer to this portion as the trunk of ${\rm Cut}(T_n)$ and denote it by
${\rm Trunk}(T_n)$. The connected components of the complement of the
trunk, ${\rm Cut}(T_n)\backslash {\rm Trunk}(T_n)$ are referred to as the
branches of ${\rm Cut}(T_n)$.

Roughly speaking, it has been shown in \cite{Be3} that upon rescaling the
graph distance of ${\rm Cut}(T_n)$ by a factor $n^{-1}\ln n$, the latter
converges to the unit interval. The precise mathematical statement involves
the notion of convergence of pointed measured metric spaces in the sense of
the Gromov-Hausdorff-Prokhorov distance.

\begin{theorem}\label{T2} Endow ${\rm Cut}(T_n)$ with the uniform
  probability measure on its leaves, and normalize the graph distance by a
  factor $n^{-1}\ln n$.
  As $n\to\infty$, the latter converges in probability in the sense of the
  pointed Gromov-Hausdorff-Prokhorov distance to the unit interval $[0,1]$
  equipped with the usual distance and the Lebesgue measure, and pointed at
  $0$.
\end{theorem}

Providing the background on the Gromov-Hausdorff-Prokhorov distance needed
to explain rigorously the meaning of Theorem \ref{T2} would probably drive
us too far away from the purpose of this survey, so we shall content
ourselves here to give an informal explanation. After the rescaling, each
edge of ${\rm Cut}(T_n)$ has length $n^{-1}\ln n$, and it follows from
\eqref{E5} that the length $n^{-1}\ln n \times L(n)$ of ${\rm Trunk}(T_n)$
converges in probability to $1$ as $n\to \infty$. Because the trunk is
merely a segment, if we equip it with the uniform probability measure on
its nodes, then we obtain a space close to the unit interval endowed with
the Lebesgue measure.  The heart of the argument of the proof in \cite{Be3}
is to observe that in turn, ${\rm Trunk}(T_n)$ is close to ${\rm Cut}(T_n)$
when $n$ is large, both in the sense of Hausdorff and in the sense of
Prokhorov. First, as ${\rm Trunk}(T_n)$ is a subset of ${\rm Cut}(T_n)$,
the Hausdorff distance between ${\rm Trunk}(T_n)$ and ${\rm Cut}(T_n)$
corresponds to the maximal depth of the branches of ${\rm Cut}(T_n)$, and
one thus have to verify that all the branches are small (recall that the
graph distance has been rescaled by a factor $n^{-1}\ln n$).  Then, one
needs to check that the uniform probability measures, respectively on the
set of leaves of ${\rm Cut}(T_n)$ and on the nodes of ${\rm Trunk}(T_n)$,
are also close to each other in the sense of the Prokhorov distance between
probability measures on a metric space. This is essentially a consequence
of the law of large numbers for the random walk defined in \eqref{ESn},
namely
$$\lim_{n\to \infty} \frac{S_n}{n\ln n}=1\qquad \hbox{in probability};$$
see \eqref{E4}. 

\subsection{Applications}
Theorem \ref{T2} enables us to specify the asymptotic behavior of the
number of cuts needed to isolate randomly chosen vertices of $T_n$. For a
given integer $\ell\geq 1$ and for each $n\geq 1$, let $U^{(n)}_1, \ldots,
U^{(n)}_{\ell}$ denote a sequence of i.i.d. uniform variables in $\{0,1,
\ldots, n\}$. We write $Y_{n,\ell}$ for the number of random cuts which are
needed to isolate $U^{(n)}_1, \ldots, U^{(n)}_{\ell}$. The following
corollary, which is taken from \cite{Be3}, is a multi-dimensional
extension of Theorem 3 of Kuba and Panholzer \cite{KP1}.

\begin{corollary} \label{C3}
As $n\to\infty$, the random vector
$$\left ( \frac{\ln n}{n} Y_{n,1}, \ldots,  \frac{\ln n}{n} Y_{n,\ell}\right)$$
converges in distribution to
$$\left( U_1, \max\{U_1, U_2\}, \ldots, \max\{U_1, \ldots, U_{\ell}\}\right),$$
where $U_1, \ldots, U_{\ell}$ are i.i.d. uniform $[0,1]$ random variables.
In particular, $\frac{\ln n}{n} Y_{n,\ell}$ converges in distribution to a
beta$(\ell,1)$ variable.
\end{corollary}

\proof Recall that $U^{(n)}_1, \ldots, U^{(n)}_{\ell}$ are $\ell$
independent uniform vertices of $T_n$. Equivalently, the singletons
$\{U^{(n)}_1\}, \ldots, \{U^{(n)}_{\ell}\}$ form a sequence of $\ell$
i.i.d. leaves of ${\rm Cut}(T_n)$ distributed according to the uniform
law. Let also $U_1, \ldots, U_l$ be a sequence of $\ell$ i.i.d. uniform
variables on $[0,1]$.  Denote by ${\mathcal R}_{n,\ell}$ the reduction of
${\rm Cut}(T_n)$ to the $\ell$ leaves $\{U^{(n)}_1\}, \ldots,
\{U^{(n)}_{\ell}\}$ and its root $\{0,1, \ldots, n\}$, i.e. ${\mathcal
  R}_{n,\ell}$ is the smallest subtree of ${\rm Cut}(T_n)$ which connects
these nodes. Similarly, write ${\mathcal R}_{\ell}$ for the reduction of
$I$ to $U_1, \ldots, U_{\ell}$ and the origin $0$. Both reduced trees are
viewed as combinatorial trees structures with edge lengths, and Theorem
\ref{T2} entails that $n^{-1}\ln n{\mathcal R}_{n,\ell}$ converges in
distribution to ${\mathcal R}_{\ell}$ as $n\to \infty$. In particular,
focusing on the lengths of those reduced trees, there is the weak
convergence
\begin{equation}
\label{E9}
\lim_{n\to \infty} \left( \frac{\ln n}{n}|{\mathcal R}_{n,1}|,\ldots,
  \frac{\ln n}{n}|{\mathcal R}_{n,\ell}|\right)  = \left( |{\mathcal
    R}_{1}|,\ldots,  |{\mathcal R}_{\ell}|\right)\qquad \hbox{in
  distribution.}
\end{equation} 
This yields our claim, as plainly $ |{\mathcal R}_{i}|= \max\{U_1, \ldots,
U_{i}\}$ for every $i=1, \ldots, \ell$. \QED

\noindent{\bf Remark.}
The nearly trivial proof of this corollary exemplifies the power of Theorem
\ref{T2}, and one might ask for convergence of the cut-tree for other tree
models. In fact, employing the work of Haas and Miermont \cite{HM}, it has been
shown in \cite{Be0} that if $T^{(c)}_n$ is a uniform Cayley tree of size
$n$, then $n^{-1/2}{\rm Cut}(T^{(c)}_n)$ converges weakly in the sense of
Gromov-Hausdorff-Prokhorov to the Brownian Continuum Random Tree (CRT), see
Aldous \cite{Al}. Since the total length of the CRT reduced to the root and
$\ell$ i.i.d leaves picked according to its mass-measure follows the
Chi$(2\ell)$-distribution, one readily obtains the statement corresponding
to Corollary \ref{C3} for uniform Cayley trees (\cite{Be0} and also, by
different means, \cite{ABBH}). Bertoin and Miermont \cite{BM} extended the
convergence of the cut-tree towards the CRT to the full family of critical
Galton-Watson trees with finite variance and conditioned to have size $n$,
in the sense of Gromov-Prokhorov. As a corollary, one obtains a
multi-dimensional extension of Janson's limit theorem \cite{Ja2}. Very
recently, Dieuleveut \cite{Di} proved the analog of \cite{BM} for the case
of Galton-Watson trees with offspring distribution belonging to the domain
of attraction of a stable law of index $\alpha\in (1,2)$.

With Corollary \ref{C3} at hand, we can also study the number $Z_{n,\ell}$
of random cuts which are needed to isolate the $\ell$ {\it last} vertices
of $T_n$, i.e. $n-\ell +1,\ldots,n$, where $\ell\geq 1$ is again a given
integer. As Kuba and Panholzer \cite{KP1} proved in their Theorem 2,
$Z_{n,\ell}$ has the same asymptotic behavior in law as $Y_{n,\ell}$.  The
following multi-dimensional version was given in \cite{Be3}, relying
on Theorem 2 of \cite{KP1}. Here we give a self-contained proof of the same
statement.

\begin{corollary} \label{C4}
As $n\to\infty$, the random vector
$$\left ( \frac{\ln n}{n} Z_{n,1}, \ldots,  \frac{\ln n}{n} Z_{n,\ell}\right)$$
converges in distribution to
$$\left( U_1, \max\{U_1, U_2\}, \ldots, \max\{U_1, \ldots, U_{\ell}\}\right),$$
where $U_1, \ldots, U_{\ell}$ are i.i.d. uniform $[0,1]$ random variables.
\end{corollary}

\proof For ease of notation, we consider only the case $\ell=1$, the
general case being similar. The random variable $Z_n=Z_{n,1}$ counts the
number of random cuts needed to isolate the vertex $n$, which is a leaf of
$T_n$. If we write $v$ for the parent of $n$ in $T_n$, then $v$ is
uniformly distributed on $\{0,1,\ldots,n-1\}$, and it follows that the
number $Y'_n$ of cuts needed to isolate $v$ has the same limit behavior in
law as $Y_{n-1,1}$. In view of Corollary \ref{C3}, it suffices therefore to
verify that
$$
\lim_{n\rightarrow\infty}\frac{\ln n}{n}\left(Y'_{n}-Z_{n}\right)=0\qquad
\hbox{in probability.}
$$
We now consider the algorithm for isolating the vertex $v$. Clearly, the
number of steps of this algorithm until the edge $e$ joining $v$ to $n$ is
removed is distributed as $Z_{n}$. In particular, we obtain a natural
coupling between $Y'_n$ and $Z_{n}$ with $Z_{n}\leq Y'_n$. Denote by
$[0;n]$ the segment of $T_n$ from the root $0$ to the leaf $n$, and write
$k$ for the outer endpoint of the first edge from $[0;n]$ which is to be
removed by the isolation algorithm. Since $|[0;n]|\sim \ln n$ in
probability (see e.g. Theorem 6.17 of \cite{Dr}), and since the isolation
algorithm chooses its edges uniformly at random, the probability that
$k$ is equal to $n$ tends to zero. Moreover, with high probability
$|[k;n]|$ will still be larger than $(\ln n)^{1/2}$, say. By conditioning on
$k$ and repeating the above argument with $[k;n]$ in place of $[0;n]$, we
see that we can concentrate on the event that before $n$ is isolated, at
least two edges different from $e$ are removed from the segment $[0;n]$.
On this event, after the second time an edge from $[0;n]$ is removed, the
vertices $v$ and $n$ lie in a tree component which can be interpreted as a
tree component of the second generation in the destruction process. As a
consequence of Proposition \ref{Pbranching}, the size of this tree
component multiplied by factor $\ln n/n$ converges to zero in
probability. Since the size of the component gives an upper bound on the
difference $Y'_n-Z_{n}$, the claim follows. \QED

As another application of the cut-tree, Theorem~\ref{T2} allows us to
determine the number of cuts $A_{n,\ell}$ which are required to {\it
  disconnect} (and not necessarily isolate) $\ell\geq 2$ vertices in $T_n$
chosen uniformly at random. For ease of description, let us assume that the
sequence of vertices $U_1^{(n)},\ldots,U_{\ell}^{(n)}$ is chosen uniformly
at random in $\{0,1,\ldots,n\}$ without replacement. Note that in the limit
$n\rightarrow\infty$, it makes no difference whether we sample with or
without replacement.

We run the algorithm for isolating the vertices
$U_1^{(n)},\ldots,U_{\ell}^{(n)}$, with the modification that we discard
emerging tree components which contain at most one of these $\ell$
vertices. We stop the algorithm when $U_1^{(n)},\ldots,U_{\ell}^{(n)}$ are
totally disconnected from each other, i.e. lie in $\ell$ different tree
components.  Write $A_{n,2}$ for the (random) number of steps of this
algorithm until for the first time, the vertices
$U_1^{(n)},\ldots,U_{\ell}^{(n)}$ do no longer belong to the same tree
component, further $A_{n,3}$ for the number of steps until for the
first time, the $\ell$ vertices are spread out over three distinct tree
components, and so on, up to $A_{n,\ell}$, the number of steps until
the $\ell$ vertices are totally disconnected. We obtain the following result.

\begin{corollary} \label{C5}
As $n\to\infty$, the random vector
$$\left ( \frac{\ln n}{n} A_{n,2}, \ldots,  \frac{\ln n}{n} A_{n,\ell}\right)$$
converges in distribution to 
$$\left( U_{(1,\ell)}, \ldots, U_{(\ell-1,\ell)}\right),$$
where $U_{(1,\ell)}\leq U_{(2,\ell)}\leq\cdots\leq U_{(\ell-1,\ell)}$ denote the first
$\ell-1$ order statistics of an i.i.d. sequence
$U_1, \ldots, U_{\ell}$ of uniform $[0,1]$ random
variables. 
\end{corollary} 
In particular, $\frac{\ln n}{n} A_{n,2}$ converges in
distribution to a beta$(1,\ell)$ random variable, and $\frac{\ln n}{n}
A_{n,\ell}$ converges in distribution to a beta$(\ell-1,2)$ law.

\proof Since the branches of ${\rm Cut}(T_n)$ are asymptotically small
compared to the trunk (see e.g. Proposition $1$ in \cite{Be3}), with
probability tending to $1$ as $n\rightarrow\infty$ the $\ell$ vertices
$U_1^{(n)},\ldots,U_{\ell}^{(n)}$ are cut from the root component one after
the other, i.e. in no stage of the disconnection algorithm, a non-root tree
component will contain more than one of the
$U_1^{(n)},\ldots,U_{\ell}^{(n)}$. On this event, writing again ${\mathcal
  R}_{n,\ell}$ for the reduction of ${\rm Cut}(T_n)$ to the $\ell$ leaves
$\{U^{(n)}_1\}, \ldots, \{U^{(n)}_{\ell}\}$ and its root
$\{0,1,\ldots,n\}$, the variable $A_{n,i+1}-1$ is given by the length of
the path in ${\mathcal R}_{n,\ell}$ from the root to the $i$th branch
point. Now, if $U_1,\ldots,U_\ell$ and ${\mathcal R}_{\ell}$ are defined as
in the proof of Corollary \ref{C3}, the distance in ${\mathcal R}_{\ell}$
from the root $0$ to the $i$th smallest among $U_1,\ldots,U_\ell$ is
distributed as $U_{(i,\ell)}$. Together with \eqref{E9}, this proves the
claim.

\noindent{\bf Remark.}
With a proof similar to that of Corollary \ref{C4}, one sees that the
statement of Corollary \ref{C5} does also hold if
$A_{n,2},\ldots,A_{n,\ell}$ are replaced by the analogous
quantities for disconnecting the $\ell$ last vertices
$n-\ell+1,\ldots,n$. On the other hand, if one is interested in
disconnecting the first $\ell$ vertices $0,\ldots,\ell -1$, and if
$B_{n,2},\ldots,B_{n,\ell}$ denote in this case the quantities
corresponding to $A_{n,2},\ldots,A_{n,\ell}$, one first
observes the trivial bound
$$
B_{n,2}\leq\cdots\leq B_{n,\ell}\leq X_{n,\ell},
$$
where $X_{n,\ell}$ is the number of steps needed to isolate
$0,1,\ldots,\ell-1$. Now, $B_{n,2}$ can be identified with the number of
steps in the algorithm for isolating the root until for the first time, an
edge connecting one of the vertices $1,\ldots,\ell-1$ to its parent is
removed. By similar means as in the proof of Lemma \ref{L1}, one readily
checks that at this time, the root component has a size of order $n^\beta$,
with $\beta$ having a beta$(\ell-1,1)$-distribution. In
particular, we see that $(X_n-B_{n,2})=o(n/\ln^2n)$ in probability, where
$X_n$ is the number of steps to isolate the root $0$. But by \eqref{E8},
also $(X_n-X_{n,\ell})=o(n/\ln^2n)$ in probability. Therefore, the
variables $B_{n,i}$ have the same limit behavior in law as $X_n$, that is
as $n\rightarrow\infty$, $\frac{\ln^2 n}{n}B_{n,i}-\ln n -\ln\ln n$,
$i=2,\ldots,\ell$, converge all to the same completely asymmetric Cauchy
variable $X$ defined by \eqref{E3}.

\subsection{Ordered destruction}
Here, we consider briefly another natural destruction procedure
of a RRT, where instead of removing edges in a uniform random order, we
remove them deterministically in their {\it natural} order. That is the
$i$th edge of $T_n$ which is removed is now the one connecting the vertex
$i$ to its parent, for $i=1,\ldots, n$.

We first point at the fact that the number of ordered edge removals which
are now needed to isolate the root (recall that we only take into account
edge removals inside the current subtree containing the
root) can be expressed as $d_n(0)=\beta_1+\cdots+\beta_n$, where
$\beta_i=1$ if the parent of vertex $i$ in $T_n$ is the root $0$, and $0$
otherwise. That is to say that $d_n(0)$ is the degree of the root. Further
the recursive construction entails the $\beta_i$ are independent variables,
such that each $\beta_i$ has the Bernoulli distribution with parameter
$1/i$. As is well-known, it then follows e.g. from Lyapunov's central
limit theorem that
$$\lim_{n\to \infty} \frac{d_n(0)-\ln n}{\sqrt {\ln n}}=\mathcal{N}(0,1)\qquad \hbox{in distribution}.$$
We refer to Kuba and Panholzer \cite{KP2} for many more results about the
degree distributions in random recursive trees.

We then turn our attention to the cut-tree described in Section
\ref{Scuttree}, which encodes the ordered destruction of $T_n$.  We write
${\rm Cut}^{\rm ord}(T_n)$ for the latter and observe that the recursive
construction of $T_n$ implies that in turn, ${\rm Cut}^{\rm ord}(T_n)$ can
also be defined by a simple recursive algorithm.  Specifically, ${\rm
  Cut}^{\rm ord}(T_1)$ is the elementary complete binary tree with two
leaves, $\{0\}$ and $\{1\}$, and root $\{0,1\}$. Once $T_n$ and hence ${\rm
  Cut}^{\rm ord}(T_n)$ have been constructed, $T_{n+1}$ is obtained by
incorporating the vertex $n+1$ and creating a new edge between $n+1$ and
its parent $U_{n+1}$, which is chosen uniformly at random in $\{0,1,\ldots,
n\}$. Note that this new edge is the last one which will be removed in the
ordered destruction of $T_{n+1}$. In terms of cut-trees, this means that
the leaf $\{U_{n+1}\}$ of ${\rm Cut}^{\rm ord}(T_n)$ should be replaced by
an internal node $\{U_{n+1},n+1\}$ to which two leaves are attached, namely
$\{U_{n+1}\}$ and $\{n+1\}$. Further, any block (internal node) $B$ of
${\rm Cut}^{\rm ord}(T_n)$ with $U_{n+1}\in B$ should be replaced by $B\cup
\{n+1\}$. The resulting complete binary tree is then distributed as
${\rm Cut}^{\rm ord}(T_{n+1})$.

If we discard labels, this recursive construction of ${\rm Cut}^{\rm
  ord}(T_n)$ corresponds precisely to the dynamics of the Markov chain on
complete binary trees described e.g. in Mahmoud \cite{Mah} for Binary
Search Trees (in short, BST). We record this observation in the following
proposition.

\begin{proposition} \label{P3} The combinatorial tree structure of ${\rm
    Cut}^{\rm ord}(T_n)$ is that of a BST with $n+1$ leaves.
\end{proposition}

BST have been intensively studied in the literature, see Drmota \cite{Dr}
and references therein, and the combination with Proposition \ref{P3}
yields a number of precise results about the number of ordered cuts which
are needed to isolate vertices in $T_n$. For instance, the so-called {\it
  saturation level} $\bar H_n$ in a BST is the minimal level of a leaf,
and can then be viewed as the smallest number of ordered cuts after which
some vertex of $T_n$ has been isolated. Similarly, the height $H_n$ is the
maximal level of a leaf, and thus corresponds to the maximal number of
ordered cuts needed to isolate a vertex in $T_n$. The asymptotic behaviors
of the height and of the saturation level of a large BST are described in
Theorem 6.47 of Drmota \cite{Dr}, in particular one has
$$\lim_{n\to \infty}\frac{\bar H_n}{\ln n} = \alpha_- \quad\hbox{and}\quad
\lim_{n\to \infty}\frac{ H_n}{\ln n} = \alpha_+$$ 
where $0<\alpha_-< \alpha_+$ are the solutions to the equation $\alpha
\ln(2\e /\alpha)=1$. In the same vein, the asymptotic results of Chauvin
{\it et al.} on the profile of large BST can be translated into sharp
estimates for the number of vertices of $T_n$ which are isolated after
exactly $k$ ordered cuts (see in particular Theorem 3.1 in \cite{CKMR}).

Finally, let us look at component sizes when edges are removed in their
natural order. Compared to uniform random edge removal, the picture is
fairly different.  Indeed, when removing an edge from $T_n$ picked
uniformly at random, the size of the subtree not containing $0$ is
distributed according to the law of $\xi$ conditioned on $\xi\leq n$. If,
in contrast, the first edge to be removed is the edge joining $1$ to its
parent $0$, then we know from \eqref{EPolya} that both originating subtrees
are of order $n$.  Since the splitting property still holds when we remove
a fixed edge, the component sizes again inherit a branching structure. In
fact, it is an immediate consequence of the definition that the structure
of the tree of component sizes corresponding to the ordered destruction on
$T_n$ agrees with the structure of $T_n$ and therefore yields the same RRT
of size $n+1$.

\section{Supercritical percolation on RRT's}
\label{Spercolation}
\subsection{Asymptotic sizes of percolation clusters}
In Section \ref{Sisoroot} it has become apparent that Bernoulli bond
percolation on $T_n$ is a tool to study the sizes of tree components which
appear in isolation algorithms. Here, we take in a certain sense the
opposite point of view and obtain results on the sizes of percolation
clusters using what we know about the sizes of tree components. Throughout
this section, we use the term cluster to designate connected components
induced by percolation, while we use the terminology tree components for
connected components arising from isolation algorithms.

More specifically, the algorithm for isolating the root can be interpreted
as a dynamical percolation process in which components that do not contain
the root are instantaneously frozen. Imagine a continuous-time version of
the algorithm, where each edge of $T_n$ is equipped with an independent
exponential clock of some parameter $\alpha$. When a clock rings, the
corresponding edge is removed if and only if it currently belongs to the
root component.  At time $t>0$, the root component can naturally be viewed
as the root cluster of a Bernoulli bond percolation on $T_n$ with parameter
$p=\exp(-\alpha t)$. Moreover, under this coupling each percolation cluster
is contained in some tree component which was generated by the isolation
process up to time $t$.  In order to discover the percolation clusters
inside a non-root tree component $T'$, the latter has to be unfrozen,
i.e. additional edges from $T'$ have to be removed. In particular, the
percolation cluster containing the root of $T'$ can again be identified as
the root component of an isolation process on $T'$, stopped at an
appropriate time.

These observations lead in \cite{Be1} to the study of the asymptotic
sizes of the largest and next largest percolation clusters of $T_n$,
when the percolation parameter $p(n)$ satisfies
\begin{equation}\label{pn}
p(n) = 1- \frac{t}{\ln n} + o(1/\ln n)\qquad \hbox{for }t>0 \hbox{ fixed.}
\end{equation}
This regime corresponds precisely to the supercritical regime, in the sense
that the root cluster is the unique giant cluster, and its complement in
$T_n$ has a size of order $n$, too. Indeed, the height $h_n$ of a vertex
$u$ picked uniformly at random in a RRT of size $n+1$ satisfies $h_n\sim\ln
n$. Since the probability that $u$ is connected to the root is given by the
first moment of $(n+1)^{-1}C_{0,n}$, where $C_{0,n}$ denotes the size of
the root cluster, one obtains
$$
\E((n+1)^{-1}C_{0,n}) = \E\left(p(n)^{h_n}\right)\sim \e^{-t}.
$$
A similar argument shows $\E((n^{-1}C_{0,n})^2)\sim \e^{-2t}$, which proves
$\lim_{n\rightarrow\infty}n^{-1}C_{0,n} = \e^{-t}$ in $L^2(\P)$.

Let us now consider the next largest clusters in the regime \eqref{pn}. We
write $C_{1,n}, C_{2,n}, \ldots$ for the sizes of the non-root percolation
clusters of $T_n$, ranked in the decreasing order.  We quote from
\cite{Be1} the following limit result.

\begin{proposition}\label{P4}
  For every fixed integer $j\geq 1$,
$$\left(\frac{\ln n}{n}C_{1,n},\ldots,\frac{\ln n}{n}C_{j,n}\right)$$
converges in distribution as $n\rightarrow\infty$ towards 
$$(x_1,\ldots,x_j),$$
where $x_1>x_2>\ldots$ denotes the sequence of the atoms of a Poisson random
measure on $(0,\infty)$ with intensity $t\e^{-t}x^{-2}\dt x$.
\end{proposition}
The intensity is better understood as the image of the intensity measure
$a^{-2}\dt a \otimes \e^{-s}\dt s$ on $(0,\infty)\times (0,t)$ by the map
$(a,s)\mapsto x=\e^{-(t-s)}a$. In fact, from our introductory remarks and
Proposition \ref{Pbranching} it should be clear that the first coordinate
of an atom $(a,s)$ stands for the asymptotic (and normalized) size of the
tree component containing the percolation cluster, while the second encodes
the time when the component was separated from the root. 

Instead of providing more details here, let us illustrate an alternative
route to prove the proposition, which was taken in~\cite{BUB} to generalize
the results to {\it scale-free random trees}. These random graphs form a
family of increasing trees indexed by a parameter $\beta \in(-1,\infty)$
that grow according to a preferential attachment algorithm, see
\cite{BA}. In the boundary case $\beta\rightarrow \infty$, one obtains a
RRT, while in the case $\beta=0$, the $i$th vertex is added to one of the
first $i-1$ vertices with probability proportional to its current
degree. In \cite{BUB}, the connection of scale-free random trees to the
genealogy of Yule processes was employed, and it should not come as a
surprise that this approach can be adapted to random recursive trees. In
fact, the case of RRT's is considerably simpler, since one has not to keep
track of the degree of vertices when edges are deleted. Let us sketch the
main changes. Denote by $T(s)$ the genealogical tree of a standard Yule
process $(\pY_r)_{r\geq 0}$ at time $s$. Similar to Section $3$
of~\cite{BUB}, we superpose Bernoulli bond percolation with parameter
$p=p(n)$ to this construction. Namely, if a new vertex is attached to the
genealogical tree, we delete the edge connecting this vertex to its parent
with probability $1-p$. We write $T^{(p)}(s)$ for the resulting
combinatorial structure at time $s$, and $T_0^{(p)}(s),
T_1^{(p)}(s),\ldots$ for the sequence of the subtrees at time $s$,
enumerated in the increasing order of their birth times, where we use the
convention that $T_j^{(p)}(s)=\emptyset$ if less than $j$ edges have been
deleted up to time $s$. In particular, $T_0^{(p)}(s)$ is the subtree
containing the root $0$, and $\sum_{i\geq 0}|T_i^{(p)}(s)| =
\pY_s$. Furthermore, if $b_i^{(p)}$ denotes the birth time of the $i$th
subtree, then the process $(T_i^{(p)}(b_i^{(p)}+s) : s\geq 0)$ is a Yule
process with birth rate $p$ per unit population size, started from a single 
particle of size $1$. By analyzing the birth times as in \cite{BUB}, one 
readily obtains the analogous statements of Section $2$ and $3$
there. This leads to another proof of Proposition \ref{P4}.

\noindent{\bf Remark.} As it is shown in the forthcoming paper~\cite{Ba},
the approach via Yule processes can be extended further to all percolation
regimes $p(n)\rightarrow 1$. Moreover, if the entire family of cluster
sizes is encoded by a tree structure similar to the tree of component
sizes, one can specify the finite-dimensional limit of this ``tree of
cluster sizes''. Details will be given in~\cite{Ba}.

\subsection{Fluctuations of the root cluster}
We finally take a closer look at the size of the root cluster $C_{0,n}$ for
supercritical percolation with parameter $$p(n) = 1-\frac{t}{\ln n}.$$
As we have already discussed, $C_{0,n}$ satisfies a law of large numbers,
but as we will point out here, $C_{0,n}$ exhibits non-Gaussian
fluctuations. This should be seen in sharp contrast to other graph models,
were asymptotic normality of the giant cluster has been established,
e.g. for the complete graph on $n$ vertices and percolation parameter
$c/n$, $c>1$ fixed (Stephanov \cite{St}, Pittel \cite{Pi}, Barraez {\it et
  al.} \cite{BBF}).

For RRT's, the fluctuations can be obtained from a recent result of
Schweinsberg \cite{Sch}. Among other things, he studied how the number of
blocks in the Bolthausen-Sznitman coalescent changes over time. The
Bolthausen-Sznitman coalescent was introduced in \cite{BS} in the context
of spin glasses, and Goldschmidt and Martin \cite{GM} discovered the
following connection to the random cutting of RRT's: Equip each edge of a
RRT of size $n$ on the vertex set $\{1,\ldots,n\}$ with an independent
standard exponential clock. If a clock rings, delete the corresponding
edge, say $e$, and the whole subtree rooted at the endpoint of $e$ most
distant from the root $1$. Furthermore, replace the label of
the vertex of $e$ which is closer to the root $1$, say $i$, by the label
set consisting of $i$ and all the vertex labels of the removed subtree.
Then the sets of labels form a partition of $\{1,\ldots,n\}$, which evolves
according to the dynamics of the Bolthausen-Sznitman coalescent started
from $n$ blocks $\{1\}, \ldots, \{n\}$ (see Proposition 2.2 of \cite{GM}
for details).

Note that in this framework, the variable $X_n$ counting the number of
steps in the algorithm for isolating the root can be interpreted as the
number of collision events which take place until there is just one block
left.

Theorem $1.7$ in \cite{Sch}, rephrased in terms of $C_{0,n}$, now
reads as follows.
\begin{theorem}{\rm (Schweinsberg \cite{Sch})}
 There is the weak convergence
 \begin{equation*}\left(n^{-1}C_{0,n}-\e^{-t}\right)\ln n -t\e^{-t}\ln\ln n
   \Longrightarrow
t\e^{-t}(X-\ln t),
\end{equation*}
where $X$ is a completely asymmetric Cauchy variable whose law is
determined by \eqref{E3}.
\end{theorem}
This statement was re-proved in \cite{Be4}, with a different approach which
does not rely on the Bolthausen-Sznitman coalescent. Instead, three
different growth phases of a RRT $T_n$ are considered, and the effect of
percolation is studied in each of these phases.  This approach makes again
use of the coupling of Iksanov and M\"ohle and the connection to Yule
processes, providing an intuitive explanation for the correction terms in
the statement.

\end{document}